\input cyracc.def
\input amssym.def
\input amssym.tex
\overfullrule 0pt
\magnification=1200
\hsize=29pc
\font\sm=cmr10 at 8pt
\font\smit= cmti10 at 8pt
\advance\vsize-14pt 
\def\today{\ifcase\month\or January\or February\or March\or April\or
May\or June\or July\or August\or September\or October\or November\or
December\fi\space\number\day, \number\year}

\def\frac#1#2{{{#1} \over {#2}}}
\def\c#1{{\cal #1}}

\def\c#1{{\cal #1}}

\font\sm=cmr10 at8pt

\centerline {\bf   An algebraic proof of Deligne's regularity criterion}
\bigskip
  \centerline{\bf by Y. Andr\'e and F.
Baldassarri}
\bigskip

\bigskip\bigskip\bigskip\bigskip\bigskip\bigskip\bigskip\bigskip\bigskip 
{\sm {\it Abstract.} Deligne's regularity criterion for an integrable 
connection $\nabla$ on a smooth complex algebraic variety $X$ says 
that $\nabla$ is regular along the irreducible divisors at infinity 
in some fixed normal compactification of $X$ if and only if the 
restriction of $\nabla$ to every smooth curve on $X$ is regular ({\it 
i. e.} has only regular singularities at infinity). The ``only if" 
part is the difficult implication. Deligne's proof is transcendental, 
and uses Hironaka's resolution of singularities. We give here an 
elementary and purely algebraic proof of this implication: it is, as far as we know, the 
first algebraic proof of Deligne's regularity criterion. 
}

\bigskip\bigskip\bigskip\bigskip\bigskip\bigskip   {\bf   Introduction.}

\medskip {\bf 0.1.} Let $X$ be a smooth connected algebraic variety 
over a field $K$ of characteristic $0$. Let ${\cal E} $ be a coherent 
module on $X$ endowed with an integrable connection $\nabla$.
When $X$ is a curve, the dichotomy between regular and irregular 
singularities (at infinity) goes back to the 19th century. The 
connection $\nabla$ is said to be {\it regular} if all singularities 
(at infinity) are regular, see Manin's classical paper [Ma]. This is 
obviously a birational notion.

\medskip In higher dimension, one may consider a normal 
compactification $\bar X$ of $X$, look at generic points $P$ of 
irreducible components of $\partial X = \bar X \setminus X$ of 
codimension one in $\bar X$, and tensor ${\cal E} $ with the discrete valuation ring ${{\cal O}_{X,P}}$. The notion of regularity of $\nabla$ {\it at 
$P$} is then the familiar one.

In order to obtain a birational notion of regularity, one is then led 
to say that $\nabla$ is {\it regular} if for any $(X', P)$ as before, 
$\nabla$ is regular at $P$.
  
  This definition (which is one of the equivalent definitions taken by Deligne 
in his Lecture Notes  [De, II.4.5]) is rather forbidding, since it 
requires to consider all divisorial valuations of $K(X)$ at the same 
time.
   Fortunately, it actually suffices to consider only one normal 
compactification $\bar X$: in fact, according to Deligne, one has the following 
characterizations of regularity:

\medskip \proclaim Theorem 0.1. [De, II.4.4, 4.6] The following are equivalent:
\medskip\noindent  $i)$ $\nabla$ is regular,
\medskip\noindent  $ii)$ for some normal compactification $\bar X$, 
$\nabla$ is regular at (the generic points of) all irreducible 
components of $\partial X$ of codimension one in $\bar X$,
\medskip\noindent  $iii)$ for any smooth curve $C$ and any locally 
closed embedding  $h: C\to X$, $h^\ast\nabla$ is regular,
\medskip\noindent  $iv)$ for any smooth $Y$ and any morphism $f: Y\to 
X$, $f^\ast\nabla$ is regular,
\medskip\noindent  $v)$ for some dominant morphism $f: Y\to X$ with 
$Y$ smooth, $f^\ast\nabla$ is regular.
\par

The difficult implication is $ii)\Rightarrow iii)$. The implication 
$iii)\Rightarrow i)$ is comparatively easy to establish ({\it cf. 
e.g.} [AB, I. 3.4.7]), and the other implications follow very easily 
from these two.
  The difficulty with $ii)\Rightarrow iii)$ arises when the closure of 
$C$ in $\bar X$ does not meet $\partial X$ transversally. A closely 
related difficulty, with $ii)\Rightarrow i)$, is to show that 
$\nabla$ remains regular at the exceptional divisor when one blows up 
a subvariety of $\partial X$.

\bigskip {\bf 0.2.}  In order to prove the difficult implication 
$ii)\Rightarrow iii)$, Deligne relied in his Lecture Notes upon a 
statement about stability of regularity under specialization. This 
statement  [De, II.1.23] is actually false, as B. Malgrange pointed 
out.

  Deligne later circulated a corrigendum, using a transcendental
argument to settle the problem (assuming, as one may, that $K=\bf 
C$), and using the existence of smooth compactifications $\bar X$ 
such that $\partial X$ is a divisor with strict normal crossings 
(Hironaka).
Over such a $\bar X$, the analytic connection $\nabla^{an}$ attached 
to $\nabla$ extends to a connection with logarithmic poles $\bar 
\nabla^{an}$ along $\partial X$, which is uniquely algebraizable. Its 
restriction to $X$ is an algebraic integrable connection 
$\nabla_{reg}$, whose analytification $\nabla_{reg}^{an}$ is 
canonically isomorphic to $\nabla^{an}$. The condition that $\nabla$ 
is regular amounts to saying that this isomorphism is algebraizable.

\bigskip {\bf 0.3.} In our book [AB], we have revisited these mostly 
well-known results about regularity, with the aim of offering purely 
algebraic, and as elementary as
possible, proofs. In particular, we claimed to give such a proof of 
the above implication $ii)\Rightarrow iii)$ ({\it cf.} [AB, I.5.4]). 
We used the fact that $\nabla$ extends to a logarithmic connection 
outside a divisor of $\partial X$ (a closed subset of codimension at 
least 2 in $\bar X$) and extended it by direct image to the whole of 
$\bar X$. We claimed that its restriction to the closure of $h(C)$ is 
a logarithmic connection, relying on a certain lemma on differentials 
with logarithmic poles (lemma 5.5).

This lemma [AB, I.5.5] is actually false\footnote{$^{(1)}$}{\sm That 
error - which we shall analyse in detail below -
questions the validity of the proof of proposition 5.4 {\smit 
loc.~cit.}. The argument sketched in the subsequent remark 5.6, as an
alternative proof of (5.4), also suffers of  the same objection.
Remark 6.5.6 of {\smit loc.~cit.}, suggested a third possible proof of
(5.4). That does not work either, since the reference to (3.4.4) is
not appropriate, because we do not have a morphism of models.} - in 
the case of a divisor with non-normal crossings -  as J. Bernstein 
pointed out to us.

\bigskip {\bf 0.4.} The aim of the present paper is twofold:

\medskip 1) to analyse in detail the nature of the error in [AB, 
I.5.5], following Bernstein (sections 3 and 4),

\medskip 2) to supply an elementary\footnote{$^{(2)}$}{\sm not using 
resolution of singularities beyond the classical case of surfaces} 
and purely algebraic proof of Deligne's regularity criterion, in the 
following slightly refined form:

    \proclaim Proposition 0.2. ({\it cf.} [AB, I.5.4]). Let $X'$ be an 
algebraic variety over a field $K$ of characteristic $0$. Let $X$ be 
a smooth open subset of $X'$, with complement ${\partial}X = 
X^{\prime} \setminus X$. We fix a closed point $P\in \partial X$ such 
that for each
irreducible component  $Z$ of ${\partial}X$ passing through $P$ whose
local ring ${\cal O}_{X^{\prime},Z}$ in $X^{\prime}$ is of dimension
$1$, that local ring is a discrete valuation ring.
  \hfill\break Let $C^{\prime} {\buildrel  h \over \longrightarrow } 
X^{\prime}$ be a
morphism from
a smooth connected $K$-curve
$C^{\prime}$ such that $h({C^{\prime}}) \not\subset {\partial}X$, and 
let $C $ be the open subset $ h^{-1}(X)$ of $C'$.
We assume there is a closed point $Q \in {C^{\prime}}$ with
$h(Q) = P$.
   \hfill\break Let  $({\cal E}, \nabla)$ be a coherent
sheaf with integrable connection on $X/K$, with generic fiber the
$\kappa(X)/K$-differential module $(E, \nabla) := ({\cal E},
\nabla)_{\eta_X}$. We assume that, for each
irreducible component  $Z$ of ${\partial}X$ passing through $P$ whose
local ring ${\cal O}_{X^{\prime},Z}$ in $X^{\prime}$ is of dimension
$1$,  $(E, \nabla)$ is regular at
the corresponding divisorial valuation $v_Z$ of $\kappa(X)$.
   \medskip Then
$({\cal E},\nabla )_{\mid C}$ is regular at $Q$.
\par
  \medskip Our proof begins with a reduction to the case of an affine 
open neighborhood $X^{\prime}$ of  $P = O$ in ${\bf
A}^2_K$ (section 1). In that situation, we proceed with a corrected 
version of the wrong lemma [AB, I.5.5], based on the study of a 
certain filtration on a logarithmic De Rham complex (section 2).

\bigskip {\it Acknowledgements.} Professor Bernstein in a
letter of September 26, 2003, indicated to us a crucial
(counter-)example on which to test a possible solution of the problem. We reproduce
a simplified version of his (counter-)example below.  His remarks
were of great help to
us. We thank him heartily.

  We are also much
indebted to Maurizio Cailotto for the key log geometry argument
in lemma 2.1.

\bigskip \noindent {\bf 1. Reduction to the two-dimensional case. }

\medskip  Our first aim is to reduce to the case when $X'$ is a
smooth
surface. This is done in the next two lemmas. This reduction serves
two purposes:
\par
- the embedded resolution of ${\partial}X $ becomes elementary,
\par
- any reflexive coherent module on such an $X'$ is locally free.
\medskip

\proclaim Lemma 1.1.  Let $X$ be a smooth
open dense subvariety of a quasi-projective irreducible $K$-variety
$X^{\prime} \subset {\bf P}^N_K$ of dimension $n > 2$, and let  $P$
be a closed point of  ${\partial}X = X^{\prime} \setminus X$.  Let
$Z_1, \dots, Z_r$ be the distinct  irreducible components  of
codimension one of  ${\partial}X$ passing through $P$.  We assume
that for any $i=1,\dots,r$ the local ring ${\cal O}_{X^{\prime},Z_i}$
is a DVR. Let
$C^{\prime} $ be a simple closed $K$-curve in $X^{\prime}$ through
$P$, ${C^{\prime}} \not\subset {\partial}X$. Let $C$ = $C^{\prime}
\cap X$. For any sufficiently large $d$, there exists an irreducible
complete intersection $Y \subset {\bf P}^N_K$ of multidegree $(d,
\dots,d)$ ($n-2$ entries) such that:
\hfill \break
$(i)$ $Y$
contains $C^{\prime}$ and cuts $X$ transversally at $\eta_C$;
\hfill
\break
$(ii)$ $Y \cap X^{\prime}$ is an irreducible surface and $Y$
cuts $X\setminus C$ transversally;
\hfill \break
$(iii)$ in a
neighborhood of $P$, $Y$ cuts each $Z_i \setminus \{ P \}$
transversally and does not cut any irreducible component of $\partial
X$ of codimension $> 1$ in $X^{\prime}$, nor the singular locus of
$\cup_i Z_i$, except in  $P$.
   \medskip
In particular, there
exists a quasi-projective irreducible neighborhood $X^{\prime}_2$ of
$P$ in $Y \cap X^{\prime}$, containing an open subset of
$C^{\prime}$, such that $X_2 :=
X^{\prime}_2 \cap X$ is smooth, the
distinct irreducible components  of codimension one of
${\partial}X_2 = X^{\prime}_2 \setminus X_2$ passing through $P$, are
precisely $Z_1 \cap X^{\prime}_2, \dots, Z_r \cap X^{\prime}_2$, and
for any $i=1,\dots,r$, the local ring ${\cal O}_{X^{\prime}_2,Z_i\cap
X^{\prime}_2}$ is a DVR.

\medskip
    \noindent
     {\bf Proof.} This is standard, but, for lack of reference, we give
some detail. Replacing $X^{\prime}$ by its closure in ${\bf
P}^N$, we may assume that $X^{\prime}$ (and $C^{\prime}$) is a closed
subvariety of ${\bf P}^N$.
Except at the very end of this proof, $n
= {\rm dim} X^{\prime}$  will be allowed to take the value $n= 2$.
Let $\pi: {\bf B} \longrightarrow {\bf P}^N$ be the blowing-up of
${\bf P}^N$ centered at $C^{\prime}$, with exceptional divisor $D$.

    We set ${\cal I}_{C^{\prime}} = {\rm Ker}
({\cal O}_{{\bf P}^N} \longrightarrow {\cal O}_{C^{\prime}}$).
    Then
$(\pi^{-1}  {\cal I}_{C^{\prime}}  ){\cal O}_{\bf B}  =  {\cal
O}(-D)$ and,   for  $d >>0$,
    $\pi^{\ast} ({\cal O}_{{\bf P}^N}(d))
\otimes
{\cal O}(-D)$ is very ample [Ha, II.7.10.$b$]: a basis of sections
defines an embedding
of  $\bf B$ into  ${\bf P}^M$. Since $\pi_{\ast}{\cal O}_{\bf B} =
{\cal O}_{{\bf P}^N}$ [Ha, proof of III.11.4], one has $\pi_{\ast}
\pi^{-1}  {\cal I}_{C^{\prime}} =    {\cal I}_{C^{\prime}}$, as
ideals of  ${\cal O}_{{\bf P}^N}$  and therefore $\pi_{\ast} {\cal
O}(-D) = {\cal I}_{C^{\prime}}$.    From the projection formula we
deduce
$$\pi_{\ast}( \pi^{\ast}({\cal O}_{{\bf P}^N}(d)) \otimes  {\cal O}(-D))
\cong {\cal O}_{{\bf P}^N}(d) \otimes {\cal I}_{C^{\prime}}
\cong {\cal I}_{C^{\prime}} (d)
\; .$$
    So, one has
    $$\eqalign {{\rm Ker} &(H^0({\bf P}^N, {\cal
O}_{{\bf P}^N}(d)) \rightarrow H^0(C^{\prime}, {\cal
O}_{C^{\prime}}(d))) = \cr
    &H^0({\bf P}^N, \, {\cal
I}_{C^{\prime}}(d)) = H^0({\bf B}, \pi^{\ast}({\cal O}_{{\bf
P}^N}(d)) \otimes {\cal O}(-D))\;  ,  \cr}
    $$
and the linear system
of hypersurfaces of degree $d $ in ${\bf P}^N$ containing $C^{\prime}$
gives rise
     to a locally closed embedding
    $${\bf P}^N \setminus C^{\prime}
\hookrightarrow {\bf P}^M = {\bf P}( H^0({\bf P}^N, {\cal I}_{C^{\prime}}(d)))
\;,$$
    with Zariski closure $\bf B$.
    The canonical bijection
between hyperplanes $\cal H$ of ${\bf P}^M$ and
    hypersurfaces $H$
of degree $d$ in ${\bf P}^N$ containing $C^{\prime}$, is such that
the intersection
${\cal H} \cap ({\bf P}^N \setminus C^{\prime})$
(in  ${\bf P}^M$)
    equals $H \setminus C^{\prime}$.  So, the
intersection of $X^{\prime} \setminus C^{\prime}$ with a general
complete intersection $Y$ of multidegree $(d,\dots,d)$ ($1 \leq s
\leq n-1$ entries) in ${\bf P}^N$ containing $\eta_C$,  is the
intersection of $X^{\prime} \setminus C^{\prime}$ with a general
linear subvariety $\cal Y$ of codimension $s$ in ${\bf P}^M$. By
[S.G.A. 4, Exp. XI, Thm. 2.1. $(i)$],
$Y$ cuts $X \setminus C$
(resp.  the smooth part of  $(\cup_i Z_i ) \setminus (C^{\prime} \cap
(\cup_i Z_i ))$) transversally  and  intersects properly  any
irreducible component of $\partial X \setminus (C^{\prime} \cap
\partial X)$.  Since $s < n$, Bertini's theorem shows that the
intersection of ${\cal Y}$ with the strict inverse image of
$X^{\prime}$ in ${\bf P}^M$ is irreducible. On the other hand, since
$\eta_C$ is a simple point of $X$, it is well-known that a general
complete intersection of $s$ hypersurfaces of degree $d$ in ${\bf
P}^N$ containing $\eta_C$,  intersects $X$ transversally at this
point.
\par
    We now apply the previous considerations for $n >2$ and
$s =n-2$ to get  properties $(i)$, $(ii)$
and $(iii)$. The last
assertion is clear.
    \medskip
    We now reconsider the situation in
the statement of the proposition. We apply lemma 1 taking as
$C^{\prime}$ the closure of the image of $h$ in $X^{\prime}$, with
reduced induced structure. The closed embedding $X^{\prime}_2
\hookrightarrow X^{\prime}$, induces, for sufficiently small open
neighborhoods $U_i$ of $\eta_{Z_i}$ in $X^{\prime}$, a closed
embedding of smooth $K$-models
$(U_i~\cap~X^{\prime}_2,U_i~\cap~Z_i~\cap X^{\prime}_2 )
\hookrightarrow (U_i ,U_i \cap Z_i)
$.  By
[AB, I.3.4.4], $({\cal E}, \nabla)_{\mid X_2}$ is regular along every
component of codimension 1 of $\partial X_2 = X^{\prime}_2 \setminus
X_2$, containing $P$. This reduces our proposition to the case where
$X^{\prime}$ is a projective surface. We next show that we may further
assume without loss of generality that  $P$ a smooth point of
$X^{\prime}$, and in fact that $X^{\prime}$ is an open neighborhood
of the origin $P = O$ in ${\bf A}^2_K$.
    \medskip \noindent

\proclaim Lemma 1.2. In the notation of lemma 1.1, let us further assume
that $X^{\prime}$ is a closed irreducible subvariety of ${\bf
P}^N_K$, of dimension 2. There exists a finite morphism
    $g:
X^{\prime} \longrightarrow {\bf P}^2_K$  whose branch locus $B
\subset {\bf P}^2_K$ does not contain the image of
$\eta_{C^{\prime}}$ and such that, for any irreducible component $W$
of $\partial X$ of dimension 1 with  $P \notin W$,  $g (P) \notin
g(W)$. \par
\noindent
    {\bf Proof.} We consider the Grassmannian  ${\bf G} :=
{\bf G}(N-3, {\bf P}^N) \cong {\bf G}(2, ({\bf P}^N)^{\vee})$ (resp.
${\bf G}(N-2, {\bf P}^N) \cong {\bf G}(1, ({\bf P}^N)^{\vee})$) of
linear subvarieties of ${\bf P}^N$ of codimension 3 (resp. 2).   Let
${\bf F} \subset  {\bf G}(N-2, {\bf P}^N) \times {\bf G}(N-3, {\bf
P}^N) $ be the incidence subvariety (locus of $({\cal L},{\alpha})$
such that ${\cal L}$ contains $\alpha$). There is a natural smooth
projective fibration $p: {\bf F} \longrightarrow {\bf G}$ with  fiber
${\bf P}^2$, where $p^{-1} (\alpha)$, for a linear subvariety
$\alpha$ of ${\bf P}^N$ of codimension 3, is the projective plane of
linear subvarieties of codimension 2 in ${\bf P}^N$, passing through
$\alpha$. The $\alpha$'s which do not intersect $X^{\prime}$, form an
open dense subset $U$ of $\bf G$. For any $\alpha \in U$, the points
of  $p^{-1} (\alpha)$ corresponding to a 2-codimensional linear
subvariety $\cal L$ of ${\bf P}^N$ which fails to intersect
$X^{\prime}$ only at smooth points and there transversally  (resp.
passing through $P$, resp. which has a non-empty intersection with
$C^{\prime}$ (resp. with $W$, for any $W$ as in the statement)), form
a Zariski closed subset $
B_{\alpha}$ (resp. $P_{\alpha}$, resp. $C_{\alpha}$ (resp.
$W_{\alpha}$)) in $p^{-1} (\alpha)$.
    All of $ B_{\alpha}$,
$C_{\alpha}$ and $W_{\alpha}$ have dimension 1, while  $P_{\alpha}$
is a point.
The curves $C_{\alpha}$ and $W_{\alpha}$ and the
closed set $P_{\alpha}$ are  irreducible. Our
problem is solved if we find $\alpha$ such that  $C_{\alpha}$ is  not
a component of $B_{\alpha}$,  and $P_{\alpha} \notin W_{\alpha}$. We
consider the  algebraic subset ${\bf H}_B$ (resp. ${\bf H}_P$, resp.
${\bf H}_C$, resp. ${\bf H}_W$) in $\bf F$ defined as the closure in
$\bf F$ of $\cup_{\alpha \in U} B_{\alpha}$ (resp.
$\cup_{\alpha \in
U} \{ P_{\alpha} \}$, resp.
$\cup_{\alpha \in U} C_{\alpha}$,
resp.
$\cup_{\alpha \in U} W_{\alpha}$). The hypersurface ${\bf H}_C$
is then irreducible, as well as the algebraic sets ${\bf H}_P$ and
${\bf H}_W$. It is clear that  ${\bf G} \cong {\bf H}_P \not\subset
{\bf H}_W$, since there are linear 2-codimensional subvarieties $\cal
L$ of ${\bf P}^N$ passing through $P$ which avoid $W$. So,  for
$\alpha$ in a dense open subset of $\bf G$, $P_{\alpha} \notin
W_{\alpha}$.   It will then suffice to prove that ${\bf H}_B$ does
not contain ${\bf H}_C$. But the curve $C^{\prime}$ is simple in
$X^{\prime}$, meaning that $C^{\prime}$ is a closed integral
$K$-subscheme of dimension 1 of $X^{\prime}$, and that
$\eta_{C^{\prime}}$ is a simple point of $X^{\prime}$.   So, let $R$
be a closed point of $C^{\prime}$, simple in
$C^{\prime}$ and in $X^{\prime}$. Then there exist linear
2-codimensional subvarieties $\cal L$ of ${\bf P}^N$, intersecting
$X^{\prime}$ transversally in $R$. This shows that ${\bf H}_B$ does
not contain ${\bf H}_C$.
\medskip

In the situation of the proposition with $X^{\prime}$ a projective
surface, we apply lemma 1.2 to the simple curve which is the closure of
the image of $h$ in $X^{\prime}$, equipped with the reduced
induced structure. We obtain an \'etale covering
$$g: V :=
X^{\prime}  \setminus (g^{-1}(B) \cup g^{-1}(g (\partial X)))
\longrightarrow
{\bf P}^2_K \setminus (B \cup g(\partial X))$$
and
push forward $({\cal E}, \nabla)$ via $g$. By [AB, I.3.4.4], the
connection $g_{\ast} (({\cal E}, \nabla)_{\mid V})$ has a regular
singularity at the generic points of the 1-dimensional components of
$B$ passing through $g(P)$, and at $g(\eta_{Z_1})$, \dots,
$g(\eta_{Z_r})$. But these are
all the generic points of the 1-dimensional irreducible components of
$B \cup g(\partial X)$ passing through $g(P)$. Therefore, the
assumptions of the proposition  are satisfied by the new choice $X :=
{\bf P}^2_K \setminus (B \cup g(\partial X))$, $X^{\prime} =  {\bf
P}^2_K $, $({\cal E}, \nabla)$ being replaced by $g_{\ast} (({\cal
E}, \nabla)_{\mid V})$, $P$ by $g(P)$, and $h$ by $g \circ h$.

Therefore we are reduced to the case when $X^{\prime}$ is a smooth
surface, and even
an affine open subset of the affine plane.

\bigskip

\noindent {\bf 2. The case of an open subset of ${\bf A}^2_K$.}
\medskip \noindent
\proclaim Lemma 2.1. Let ${X^{\prime}}= {\rm Spec} A$, be an
affine open $K$-neighborhood of the origin $O$ in ${\bf A}^2_K$ and
set $U := {X^{\prime}} \setminus \{ O \}$. Let $g_1, \dots, g_r $ be
non-associated  irreducible elements of $A$,  $Z_i = V(g_i)$, for $i
=1, \dots, r$, and $Z:= \cup_i Z_i$. We assume that $U \cap Z$ is a
normal crossing divisor in $U$.
\hfill\break Let
$C^{\prime} {\buildrel  h \over \longrightarrow } X'$ be a
    morphism from
a smooth connected $K$-curve
$C^{\prime}$
and assume there is a closed point $Q \in {C^{\prime}}$ with
$h(Q) = O$, while, for $C := C^{\prime} \setminus \{ Q \}$, $h(C)
\subset X$.
\hfill\break If $({\cal E}, \nabla)$ is an
integrable connection on ${X} ={X^{\prime}} \setminus Z$, regular at
$\eta_{Z_1}$,
\dots, $\eta_{Z_r}$, then $({\cal E},\nabla )_{\mid C}$ is regular at $Q$.
\par \noindent

{\bf Proof.}

We are presently working on this section!

\bigskip  \noindent

{\bf 3. Example (after J. Bernstein).}
\medskip
We consider the situation of lemma 2.1. Let $(u,v)$ be affine
coordinates on  $X^{\prime} = {\bf A}^2$, $Z := V(uv(u+v))$. Let $\pi
: \tilde{X} \longrightarrow X^{\prime} $ be the
blowing-up of $X^{\prime} $ centered at $O$. We use a local coordinate $t =
\frac{v}{u}$ on $\tilde{X}$.

\medskip \noindent {\bf 3.1.}  $\omega^1_{(X^{\prime} ,Z)} \subsetneqq
j_{\ast}\Omega^1_U(\log Z \cap U) $
\medskip
     Consider the
$1$-form
     $$ \quad \alpha = \frac{du}{u} + \frac{dv}{v} -2
\frac{d(u+v)}{u+v} \in \Gamma (X, j_{\ast}\Omega^1_U(\log Z \cap U))
\leqno (3.1.1)$$
   An easy computation shows that
   $$ \quad    uv(u+v)\alpha = (u - v)(u dv - v du) \leqno (3.1.2)$$
   Now set $\omega = \frac{\alpha}{u-v}$. To show that $\omega \in
\Gamma (X, j_{\ast}\Omega^1_U(\log Z \cap U))$, we use the criterion
of Deligne [De, I.2.2.1].
Consider the behavior of the form $\omega$ on $U$. Formula (3.1.2)
    implies that $uv(u+v) \omega = u dv - v du$. Hence
    the form $\omega$ is regular on $X \setminus Z$. In a neighborhood
of the variety $Z \bigcap U$ the form $\omega$ is proportional to the
form $\alpha$ which by construction has simple poles along $Z \bigcap
U$.  This
shows that the form $\omega$ has simple poles along $Z \cap U$ on
$U$.  On the other hand, $d \omega = - \frac{du \wedge dv}{uv(u+v)}$
also has only simple poles along $Z \cap U$. We conclude that
$\omega$ is a section of $j_{\ast} \Omega^1_U(\log Z \cap U)$. Notice
that $\omega$ is not a section of $\omega^1_{(X^{\prime} ,Z)}$, since
$\pi^{\ast}\omega =
\frac{dt}{ut(1+t)} $, so that $d\pi^{\ast}\omega = - \frac{du \wedge
dt}{u^2t(1+t)} $ has a double pole along the exceptional divisor, and
is not therefore a section of $\Omega^1_{\tilde{X}}(\log
\tilde{Z})$.
A more precise computation shows that
$\Gamma (U,
\Omega^1_U(\log Z \cap U))$ is the free $K[u,v]$-module generated by
$\omega$ and $\frac{du}{u}$ (or by $\omega$ and $\frac{dv}{v}$).

On
the other hand, the closed exact immersion of the log scheme
$(X^{\prime} ,Z)$ in the smooth (over $K$) log scheme
$({\bf A}^3,
V(xyz))$, $(u,v) \longmapsto (u,v,u+v)$, shows  [IL1, II.2.5.$(a)$]
that
$\omega^1_{(X^{\prime} ,Z)}$ is the $K[u,v]$-module generated by
${\rm dlog}(u)$, ${\rm dlog}(v)$, ${\rm dlog}(u+v)$, subject to the
single relation
$$u\,{\rm dlog}(u) + v\,{\rm dlog}(v)  = (u+v) \,{\rm
dlog}(u+v).
$$

     \medskip
\noindent  {\bf 3.2. A counterexample to [AB, I.5.5].}
\medskip

Let now $C^{\prime} $ be the curve described in parametric form by $u 
= s, v = \nu s^2$, for generic
$\nu \in K$.

  The restriction of  the form $\omega$ (which is a section of 
$j_{\ast} \Omega^1_U(\log Z \cap U)$) to the curve $C^{\prime} $ 
equals $
\frac{ds}{s^2} (1 + O(s))$  and thus is {\it not of logarithmic type} at $s
=0$.

This is a counter-example to [AB, I.5.5].

\medskip

A consequence of the previous construction is that the restriction to 
$C^{\prime} $ of the direct image by $j$ of an integrable connection 
with logarithmic poles along $Z\cap  U$ is not necessarily of 
logarithmic type at $s=0$.

Explicitly, let $\c E = \c O^2_U$ be equipped with the
logarithmic connection $d + G$, where
$$G = \pmatrix{0&0\cr
\omega&\frac{du}{2u} + \frac{dv}{2v}} \; .
$$
Since $dG + G \wedge
G = 0$, the connection is integrable. What we said before about $\omega$
shows however that the connection given on $j_{\ast}\c E = \c
O^2_{X^{\prime}}$ by the same formula is not logarithmic.  The matrix 
of the pulled-back connection on $\c
O^2_{C^{\prime}}$ has a {\it double pole at $s=0$}.
\par To show that this
connection is nevertheless {\it regular}, we follow the proof of lemma 3 above.
The
pulled-back connection on $\tilde{\c E} = \c O_{\tilde{X}}^2$ is
locally expressed as $d+H$, where
$$H = \pmatrix{0&0\cr
\frac{dt}{ut(1+t)}&\frac{du}{u} + \frac{dt}{2t}} \; .
$$
The
restriction of this connection to the subsheaf $\c F = u \c
O_{\tilde{X}} \oplus \c O_{\tilde{X}}$ of $\tilde{\c E}$, meaning, in
matrix notation,
$$ \c F = \{ \pmatrix{uf \cr g \cr} \mid f,g \in \c
O_{\tilde{X}} \}   ,$$
is therefore logarithmic. We now pull-back
$(\c F, d +H)$ via the strict inverse image of $C^{\prime}$, namely
$\tilde{h} : C^{\prime} \hookrightarrow \tilde{X}$,  $s \longmapsto
(u,t) = (s, \nu s)$, and we get that
$$ \tilde{h}^{\ast} \c F = \{
\pmatrix{sf \cr g \cr} \mid f,g \in \c O_{C^{\prime}} \}  \; ,$$
is
stable under $d + H(s)$, where
$$H(s) = \pmatrix{0&0\cr
\frac{ds}{s^2(1+\nu s)}&\frac{3ds}{2s}} \; .
$$

\bigskip \bigskip\vfill \eject
\centerline{\bf References}
\bigskip
\item{[AB]} Andr\'e Y., Baldassarri F. : {\it ``De Rham Cohomology of
Differential Modules on Algebraic Varieties"},
Progress in Mathematics, Vol. 189, Birkh\"auser Verlag 2001.
\medskip
\item{[De]} Deligne P.: {\it ``Equations Diff\'erentielles \`a
Points Singuliers R\'eguliers"}, Lecture Notes in
Math.  163, Springer-Verlag 1970; Erratum, dated April 1971.
\medskip
\item{[Ha]} Hartshorne R. : {\it ``Algebraic Geometry"},
Springer Graduate Texts in Mathematics 52, Springer Verlag 1977.
\medskip
\item {[IL1]} Illusie L.: ``G\'eom\'etrie logarithmique, Expos\'e I-II", notes
of a
course held during the $p$-adic semester at IHP Paris in Spring
1997.
   \medskip
\item {[IL2]} Illusie L.: ``An overview of the work of K. Fujiwara, K.
Kato, and C. Nakayama on logarithmic \'etale cohomology",
in {\it Cohomologies $p$-adiques et Applications Arithm\'etiques"},
Ast\'erisque 279 (2002),  271-322.
\medskip
\item{[KK]} Kato K.: ``Logarithmic structures of Fontaine-Illusie",
in {\it ``Algebraic Analysis, Geometry and Number Theory"},
The Johns Hopkins Univ. Press (1989), 191-224.
\medskip
\item {[Ma]} Manin Yu.: ``Moduli fuchsiani", Annali Scuola Norm. Sup. 
Pisa III 19 (1965), 113-126.
\medskip
\item {[Se]} Serre J.P.: ``Prolongement de faisceaux analytiques coh\'erents",
Ann. Inst. Fourier Grenoble 16,1 (1966), 363-374.

\bigskip
\settabs 2 \columns

\+{\it Yves Andr\'e} & {\it Francesco Baldassarri} & \cr
\+
D.M.A.& Dipartimento di Matematica& \cr
\+ \' Ecole Normale Sup\'erieure & Universit\`a di Padova & \cr
\+ 45 rue d'Ulm &
Via Belzoni 7& \cr
\+
F-75230 Paris - France
&I-35131 Padova - Italy & \cr
\+ e-mail: andre@dma.ens.fr &
e-mail: baldassa@math.unipd.it& \cr
\vfill

\end